\documentclass[12pt]{article}
\usepackage{amsfonts}
\usepackage{amssymb}

\input{tcilatex}
\begin{document}

\begin{center}
\textbf{Effros, Baire, Steinhaus and Non-Separability}\\[0pt]
\bigskip

\textbf{By A. J. Ostaszewski}

\bigskip
\end{center}

\noindent \textbf{Abstract. }We give a short proof of an improved version of
the Effros Open Mapping Principle via a shift-compactness theorem (also with
a short proof), involving `sequential analysis' rather than separability,
deducing it from the Baire property in a general Baire-space setting (rather
than under topological completeness). It is applicable to
absolutely-analytic normed groups (which include complete metrizable
topological groups), and via a Steinhaus-type Sum-set Theorem (also a
consequence of the shift-compactness theorem) includes the classical Open
Mapping Theorem (separable or otherwise).

\bigskip

\textbf{Keywords:} Open Mapping Theorem, absolutely analytic sets, base-$%
\sigma $-discrete maps, demi-open maps, Baire spaces, Baire property,
group-action shift-compactness.

\textbf{Classification Numbers}: 26A03; 04A15; 02K20.

\section{Introduction}

\noindent We generalize a classic theorem of Effros \cite{Eff} beyond its
usual separable context. Viewed, despite the separability, as a group-action
counterpart of the Open Mapping Theorem OMT (that a surjective continuous
linear map between Fr\'{e}chet\textit{\ }spaces is open -- cf. \cite{Rud}),
it has come to be called the \textit{Open Mapping Principle} -- see \cite[\S %
1]{Anc}. Our `non-separable' approach is motivated by a sequential property
related to the Steinhaus-type Sum-set Theorem (that $0$ is an interior point
of $A-A,$ for non-meagre $A$ with BP, the Baire property --\cite{Oxt}, \cite%
{Pic}), because of the following argument (which goes back to Pettis \cite%
{Pet}).

Consider $L:E\rightarrow F,\ $a linear, continuous surjection between Fr\'{e}%
chet spaces, and $U$ a neighbourhood (nhd) of the origin. Choose $A$ an 
\textit{open} nhd of the origin with $A-A\subseteq U$; as $L(A)$ is
non-meagre (since $\{nL(A):n\in \mathbb{N}\}$ covers $F$) and has BP (see
Proposition 2 in \S 2.3), $L(A)-L(A)$ is a nhd of the origin by the Sum-set
Theorem. But of course%
\[
L(U)\supseteq L(A)-L(A), 
\]%
so $L(U)$ is a nhd of the origin. So $L$ is an open mapping.\footnote{%
This proof is presumably well-known -- so simple and similar to that for the
automatic continuity of homomorphisms -- but we have no textbook reference;
cf. \cite[Cor. 1.5]{KalPR}.}

\bigskip

Throughout this paper, without further comment, all spaces considered will
be metrizable, but not necessarily separable. We recall the
Birkhoff-Kakutani theorem (cf. \cite[\S II.8.3]{HewR}), that a metrizable
group $G$ with neutral element $e_{G}$ has a right-invariant metric $%
d_{R}^{G}$. Passage to $||g||:=d_{R}^{G}(g,e_{G})$ yields a (group) norm
(invariant under inversion, satisfying the triangle inequality, cf. \cite%
{ArhT}, \cite{BinO1} -- see \S 4.1), which justifies calling these \textit{%
normed groups}; any Fr\'{e}chet space qua additive group, equipped with an
F-norm (\cite[Ch. 1 \S 2]{KalPR}), is a natural example (cf. $Auth$ in \S %
2.2). Recall that a \textit{Baire space} is one in which Baire's theorem
holds -- see \cite{AaL}. Below we need the following.

\bigskip

\noindent \textbf{Definitions 1} (cf. \cite{Pet}). For $G$ a metrizable
group, say that $\varphi :G\times X\rightarrow X$ is a \textit{Nikodym }group%
\textit{\ }action (or that it has the Nikodym property -- cf. \cite{BinO4})
if for every non-empty open neighbourhood $U\ $of $e_{G}$ and every $x\in X$
the set $Ux=\varphi _{x}(U):=\varphi (x,U)$ contains a non-meagre \textit{%
Baire set}. (Here Baire set, as opposed to Baire space as above, means `set
with the Baire property'.)

\noindent\ \textbf{2. }$A^{q}$ denotes the \textit{quasi-interior} of $A$ --
the largest open set $U$ with $U\backslash A$ meagre (cf. \cite[\S 4]{Ost1}%
); other terms (`analytic', `base-$\sigma $-discrete', `group action') are
recalled later.

\bigskip

Concerning when the above property holds see \S 2.3. Our main results are
Theorems S and E below, with Corollaries in \S 2.3 including OMT; see below
for commentary.

\bigskip

\textbf{Theorem S (Shift-compactness Theorem).} \textit{For }$T$\textit{\ a
Baire non-meagre subset of a metric space }$X$\textit{\ and }$G$\textit{\ a
group, Baire under a right-invariant metric, and with separately continuous
and transitive Nikodym action on }$X$\textit{:}

\textit{for every convergent sequence }$x_{n}$\textit{\ with limit }$x$%
\textit{\ and any Baire non-meagre }$A\subseteq G$ \textit{with }$e_{G}\in
A^{q}$ \textit{and} $A^{q}x\cap T^{q}\neq \emptyset ,$\textit{\ there are }$%
\alpha \in A$\textit{\ and an integer }$N$\textit{\ such that }$\alpha x\in
T $\textit{\ and}%
\[
\{\alpha (x_{n}):n>N\}\subseteq T. 
\]%
\textit{In particular, this is so if }$G$ \textit{is analytic and all
point-evaluation maps }$\varphi _{x}$ \textit{are base-}$\sigma $-\textit{%
discrete.}

\bigskip

This theorem has wide-ranging consequences, including Steinhaus' Sum-set
Theorem -- see the survey article \cite{Ost4}, and the recent \cite{BinO3};
cf. \S 4.2-4.

\bigskip

\textbf{Theorem E (Effros Theorem -- Baire version).} \textit{If}\newline
(i)\textit{\ the normed group }$G\ $\textit{has separately continuous and\
transitive Nikodym action on }$X$\textit{;\newline
}(ii)\textit{\ }$G$ \textit{is Baire under the norm topology and }$X$\textit{%
\ is non-meagre\newline
-- then for any open neighbourhood }$U$\textit{\ of }$e_{G}$\textit{\ and
any }$x\in X$\textit{\ the set }$Ux:=\{u(x):u\in U\}$ \textit{is a
neighbourhood of }$x$\textit{, so that in particular the point-evaluation
maps }$g\rightarrow g(x)$ \textit{are open for each }$x$\textit{. That is,
the action of }$G$\textit{\ is micro-transitive. }

\textit{In particular, this holds if }$G$\textit{\ is analytic and Baire,
and all point-evaluation maps }$\varphi _{x}$ \textit{are base-}$\sigma $%
\textit{-discrete.}

\bigskip

By Proposition B2 (\S 2.3) $X$, being non-meagre here, is also a Baire space.

The classical counterpart of Theorem E has $G$ a Polish group; van Mill's
version \cite{vMil1} requires the group $G$ to be analytic (i.e. the
continuous image of some Polish space, cf. \cite{JayR}, \cite{Kec2}). The
Baire version above improves the version given in \cite{Ost3}, where the
group is almost complete. (The two cited sources taken together cover the
literature.)

A result due to Loy \cite{Loy} and to Hoffmann-J\o rgensen \cite[Th. 2.3.6
p. 355]{HofJ} asserts that a Baire, separable, analytic \textit{topological} 
\textit{group} is Polish (as a consequence of an analytic group being
metrizable -- for which see again \cite[Th. 2.3.6]{HofJ}), so in the
analytic separable case Theorem E reduces to its classical version.

Unlike the proof of the Effros Theorem attributed to Becker in \cite[Th. 3.1]%
{Kec1}, the one offered here does not employ the Kuratowski-Ulam Theorem
(the Category version of the Fubini Theorem), a result known to fail beyond
the separable context (as shown in \cite{Pol}, cf. \cite{vMilP}, but see 
\cite{FreNR}).

For recent work on the circumstances when Theorems E and S are equivalent,
see [BinO4]. For further commentary (connections between convexity and the
Baire property, relation to van Mill's separation property in \cite{vMil2},
certain specializations) see \S 4.

\section{Analyticity, micro-action, shift-compactness}

We recall some definitions from general topology, before turning to ones
that are group-related. We refer to \cite{Eng} for general topological usage
(but prefer `meagre' to `of first category'); see also \S 4.5,6.

\subsection{Analyticity}

We say that a subspace $S\ $of a metric space $X$ has a \textit{Souslin}-$%
\mathcal{H}$ \textit{representation} if there is a \textit{determining system%
} $\langle H(i|n)\rangle :=\langle H(i|n):i\in \mathbb{N}^{\mathbb{N}%
}\rangle $ of sets in $\mathcal{H}$ with (\cite{Rog}, \cite{Han1}, \cite%
{Han2}) 
\[
S=\bigcup\nolimits_{i\in I}\bigcap\nolimits_{n\in \mathbb{N}}H(i|n),\text{ }%
(I:=\mathbb{N}^{\mathbb{N}},\text{ }i|n:=(i_{1},...,i_{n})).
\]%
A topological space is an (absolutely) \textit{analytic} space if it is
embeddable as a Souslin-$\mathcal{F}$ set in its own metric completion (with 
$\mathcal{F}$ the closed sets); in particular, in a complete metric space $%
\mathcal{G}_{\delta }$-subsets (being $\mathcal{F}_{\sigma \delta }$) are
analytic. For more recent generalizations see e.g. \cite{NamP}. According to
Nikodym's theorem, if $\mathcal{H}$ above comprises Baire sets, then also $S$
is Baire (the Baire property is preserved by the Souslin operation): so
analytic subspaces are Baire sets. For background -- see \cite{Kec2} Th.
21.6 (the Lusin-Sierpi\'{n}ski Theorem) and the closely related Cor. 29.14
(Nikodym Theorem), cf. the treatment in \cite{Kur} Cor. 1 p. 482, or \cite%
{JayR} pp. 42-43. For the extended Souslin operation of non-separable
descriptive theory see also \cite{Ost2}. This motivates our interest in
analyticity as a carrier of the Baire property, especially as continuous
images of separable analytic sets are separable, hence Baire.

However, the continuous image of an analytic space is not in general
analytic -- for an example of failure see \cite{Han3} Ex. 3.12. But this
does happen when, additionally, the continuous map is base-$\sigma $%
-discrete, as defined below (\textit{Hansell's Theorem}, \cite{Han3} Cor.
4.2)\textbf{. }This technical condition is the standard assumption for
preservation of analyticity and holds automatically in the separable realm.
Special cases include \textit{closed surjective} maps and \textit{%
open-to-analytic injective} maps (taking open sets to analytic sets). To
define the key concept just mentioned, recall that for an (indexed) family $%
\mathcal{B}:=\{B_{t}:t\in T\}$:

\noindent (i) $\mathcal{B}$ is \textit{index-discrete} in the space $X$ (or
just \textit{discrete} when the index set $T$ is understood) if every point
in $X$ has a nhd meeting the sets $B_{t}$ for at most one $t\in T,$

\noindent (ii) $\mathcal{B}$ is $\sigma $-\textit{discrete} if $\mathcal{B}%
=\bigcup\nolimits_{n}\mathcal{B}_{n}$ where each set $\mathcal{B}_{n}$ is
discrete as in (i), and

\noindent (iii) $\mathcal{B}$ is a \textit{base} \textit{for} $\mathcal{A}$
if every member of $\mathcal{A}$ is the union of a subfamily of $\mathcal{B}$%
. For $\mathcal{T}$ a topology (the family of all open sets) with $\mathcal{B%
}\subseteq \mathcal{T}$ a base for $\mathcal{T}$, this reduces to $\mathcal{B%
}$ being simply a (topological) \textit{base}.

\bigskip

\textbf{Definitions.} 1.~(\cite{Mic1}, Def. 2.1) Call $f:X\rightarrow Y$ 
\textit{base-}$\sigma $\textit{-discrete} (or \textit{co-}$\sigma $\textit{%
-discrete, }\cite[\S 3]{Han3}) if the image under $f$ of any discrete family
in $X$ has a $\sigma $-discrete base in $Y.$\newline
2 (\cite[\S 2]{Han3}, cf. \cite{KanP}). An indexed family $\mathcal{A}%
:=\{A_{t}:t\in T\}$ is $\sigma $-\textit{discretely decomposable} if there
are discrete families $\mathcal{A}_{n}:=\{A_{tn}:t\in T\}$ such that $%
A_{t}=\bigcup\nolimits_{n}A_{tn}$ for each $t.$\newline
3 (\cite{Mic1}, Def. 3.3). Call $f:X\rightarrow Y$ \textit{index-}$\sigma $%
\textit{-discrete} if the image under $f$ of any discrete family $\mathcal{E}
$ in $X$ is $\sigma $-discretely decomposable in $Y.$ (Note $f(\mathcal{E)}$
is regarded as indexed by $\mathcal{E}$, so could be discrete without being
index-discrete.)

\subsection{Action, micro-action, shift-compactness}

Recall that a normed group $G\ $\textit{acts continuously} on $X$ if there
is a continuous mapping $\varphi :G\times X\rightarrow X$ such that $\varphi
(e_{G},x)=x$ and $\varphi (gh,x)=\varphi (g,\varphi (h,x))$ ($x\in X,g,h\in
G).$The action $\varphi $ is \textit{separately continuous} if $g:x\mapsto
\varphi (g,x)$ is continuous for each $g,$ and $\varphi _{x}:g\mapsto
\varphi (g,x)$ is continuous for each $x$; in such circumstances:

\noindent (i) the elements $g\in G$ yield autohomeomorphisms of $X$ via $%
g:x\mapsto g(x):=\varphi (g,x)$ (as $g^{-1}$ is continuous), and

\noindent (ii) point-evaluation of these homeomorphisms, $\varphi
_{x}(g)=g(x),$ is continuous.

\noindent In certain situations joint continuity of action is implied by
separate continuity (see \cite{Bou} and literature cited in \cite{Ost2}).

The action is \textit{transitive} if for any $x,y$ in $X$ there is $g\in G$
such that $g(x)=y.$ For later purposes (\S 2.3 and 3), say that the action
of $G$ on $X$ is \textit{weakly micro-transitive} if for $x\in X$ and each
nhd $A$ of $e_{G}$ the set 
\[
\text{cl}(Ax)=\text{cl}\{ax:a\in A\} 
\]%
has $x$ as an interior point (in $X)$. The action is \textit{micro-transitive%
} (`transitive in the small' -- for details see \cite{vMil1}) if for $x\in X$
and each nhd $A$ of $e_{G}$ the set 
\[
Ax=\{ax:a\in A\} 
\]%
is a nhd of $x.$ This (norm) property implies that $Ux$ is open for $U$ open
in $G$ (i.e. that here each $\varphi _{x}$ is an open mapping). We refer to $%
Ax$ as an $x$ \textit{orbit} (the $A$-orbit of $x$). The following group
action connects the Open Mapping Theorem to the present context.

\bigskip

\textbf{Example (Induced homomorphic action}). A surjective, continuous
homomorphism $\lambda :G\rightarrow H$ between normed groups induces a
transitive action of $G$ on $H$ via $\varphi ^{\lambda }(g,h):=\lambda (g)h$
( cf. \cite{Ost2} Th. 5.1), specializing for $G,H$ Fr\'{e}chet spaces
(regarded as normed, additive groups) and $\lambda =L:G\rightarrow H$ linear
(Ancel \cite{Anc} and van Mill \cite{vMil1}) to%
\[
\varphi ^{L}(a,b):=L(a)+b. 
\]%
Of course for Fr\'{e}chet spaces, by the Open Mapping Theorem itself, $%
\varphi ^{L}$ has the Nikodym property.

\bigskip

\textbf{Definitions. }1. $Auth(X)$ denotes the autohomeomorphisms of a
metric space $(X,d^{X})$; this is a group under composition. $\mathcal{H}(X)$
comprises those $h\in Auth(X)$ of bounded norm:%
\[
||h||:=\sup\nolimits_{x\in X}d^{X}(h(x),x)<\infty . 
\]

\noindent 2. For a normed group $G$ acting on $X,$ say that $X$ has the 
\textit{crimping property} (property C for short) w.r.t. $G$ if, for each $%
x\in X$ and each sequence $\{x_{n}\}\rightarrow x,$ there exists in $G$ a
sequence $\{g_{n}\}\rightarrow e_{G}$ with $g_{n}(x)=x_{n}.$ (This and a
variant occurs in \cite[Ch. III; Th.4]{Ban}; and \cite{ChCh}; for the term
see \cite{BinO2}.)

For a subgroup $\mathcal{G\subseteq H}(X)$, say that $X$ has the \textit{%
crimping property} w.r.t. $\mathcal{G}$ if $X$ has the crimping property
w.r.t. to the natural action $(g,x)\rightarrow g(x)$ from $\mathcal{G}\times
X\rightarrow X.$ (This action is continuous relative to the left or right
norm topology on $\mathcal{G}$ -- cf. \cite{Dug} XII.8.3, p. 271.)

\noindent 3. As a matter of convenience, say that the \textit{Effros property%
} (or \textit{property E}) holds for the group $G\ $acting on $X$ if the
action is micro-transitive, as above.

\noindent 4. For a subgroup $\mathcal{G}\subseteq Auth(X)$ say that $X$ is $%
\mathcal{G}$\textit{-shift-compact} (or, shift-compact under $\mathcal{G}$)
if for any convergent sequence $x_{n}\rightarrow x_{0},$ any open subset $U$
in $X$ and any Baire set $T$ co-meagre in $U,$ there is $g\in \mathcal{G}$
with $g(x_{n})\in T\cap U$ along a subsequence. Call the space \textit{%
shift-compact }if it is $\mathcal{H}(X)$-shift-compact (cf. \cite{MilO}, 
\cite{Ost5}).

In such a space, any Baire non-meagre set is locally co-meagre (co-meagre on
open sets) in view of Prop. B2 below.

We shall prove in \S\ 3.1 equivalence between the Effros and Crimping
properties:

\bigskip

\textbf{Theorem EC.} \textit{The Effros property holds for a group }$G$%
\textit{\ acting on }$X$\textit{\ iff }$X$\textit{\ has the Crimping
property w.r.t. }$G.$

\bigskip

We now clarify the role of shift-compactness.

\bigskip

\textbf{Proposition B1.} \textit{For any subgroup }$\mathcal{G}\subseteq 
\mathcal{H}(X),$\textit{\ if }$X$\textit{\ is }$\mathcal{G}$\textit{%
-shift-compact, then }$X$\textit{\ is a Baire space.}

\bigskip

\textbf{Proof.} We argue as in \cite{vMil2} Prop 3.1 (1). Suppose otherwise;
then $X$ contains a non-empty meagre open set. By Banach's Category Theorem
(or localization principle, for which see \cite{JayR} p. 42, or \cite{Kel}
Th. 6.35), the union of all such sets is a largest open meagre set $M,$ and
is non-empty. Thus $X\backslash M$ is a co-meagre Baire set. For any $x\in M$
the constant sequence $x_{n}\equiv x$ is convergent and, since $X\backslash
M $ is co-meagre in $X,$ there is $g\in G$ with $g(x)\in X\backslash M$.
But, as $g$ is a homeomorphism, $g(M)$ is a non-empty open meagre set, so is
contained in $M,$ implying $g(x)\in M,$ a contradiction. $\square $

\bigskip

A similar argument gives the following and clarifies an assumption in
Theorem E.

\bigskip

\textbf{Proposition B2 }(cf. \cite{vMil2}; \cite[Prop. 2.2.3]{HofJ}).\textbf{%
\ }\textit{If }$X$\textit{\ is non-meagre and }$G$\textit{\ acts
transitively on }$X,$\textit{\ then }$X$\textit{\ is a Baire space.}

\bigskip

\textbf{Proof.} As above, refer again to $M,$ the union of all meagre open
sets, which, being meagre, has non-empty complement. For $x_{0}$ in this
complement and any non-empty open $U$ pick $u\in U$ and $g\in G$ such that $%
g(x_{0})=u.$ Now, as $g$ is continuous, $g^{-1}(U)$ is a nhd of $x_{0},$ so
is non-meagre, since every nhd of $x_{0}$ is non-meagre. But $g$ is a
homeomorphism, so $U=g(g^{-1}(U)$ is non-meagre. So $X$ is Baire, as every
non-empty open set is non-meagre. $\square $

\subsection{Nikodym actions}

The following result generalizes one that, for separable groups $G$, is
usually a first step in proving the weakly micro-transitive variant of the
classical Effros Theorem (cf. Ancel \cite{Anc} Lemma 3, \cite{Ost3} Th. 2).
Indeed, one may think of it as giving a form of `very weak
micro-transitivity'.

\bigskip

\textbf{Proposition 1.} \textit{If }$G$\textit{\ is a normed group, acting
transitively on a non-meagre space }$X$\textit{\ with each point evaluation
map }$\varphi _{x}:g\mapsto g(x)$\textit{\ base-}$\sigma $\textit{-discrete
relative -- then for each non-empty open }$U$\textit{\ in }$G$ \textit{and
each }$x\in X$\textit{\ the set }$Ux$\textit{\ is non-meagre in }$X$\textit{%
. }

\textit{In particular, if }$G$\textit{\ is analytic, then }$G$\textit{\ is a
Nikodym action.}

\bigskip

\textbf{Proof.} We first work in the right norm topology, i.e. derived from
the assumed right-invariant metric $d_{R}^{G}(s,t)=||st^{-1}||$. Suppose
that $u\in U,$ and so without loss of generality assume that $%
U=B_{\varepsilon }(u)=B_{\varepsilon }(e_{G})u$ (open balls of radius some $%
\varepsilon >0);$ then put $y:=ux$ and $W=B_{\varepsilon }(e_{G}).$ Then $%
Ux=Wy.$ Next work in the left norm topology, derived from $%
d_{L}^{G}(s,t)=||s^{-1}t||=d_{R}^{G}(s^{-1},t^{-1})$ (for which $%
W=B_{\varepsilon }(e_{G})$ is still a nhd of $e_{G})$. As each set $hW$ for $%
h\in G$ is now open (since now the left shift $g\rightarrow hg$ is a
homeomorphism), the open family $\mathcal{W}=\{gW:g\in G\}$ covers $G.$ As $%
G $ is metrizable (and so has a $\sigma $-discrete base), the cover $%
\mathcal{W}$ has a $\sigma $-discrete refinement, say $\mathcal{V}%
=\bigcup\nolimits_{n\in \mathbb{N}}\mathcal{V}_{n}$, with each $\mathcal{V}%
_{n}$ discrete. Put $X_{n}:=\bigcup \{Vy:V\in \mathcal{V}_{n}\};$ then $%
X=\bigcup\nolimits_{n\in \mathbb{N}}X_{n},$ as $X=Gy,$ and so $X_{n}$ is
non-meagre for some $n,$ for $n=N$ say. Since $\varphi _{y}$ is base-$\sigma 
$-discrete, $\{Vy:V\in \mathcal{V}_{N}\}$ has a $\sigma $-discrete base, say 
$\mathcal{B}=\bigcup\nolimits_{m\in \mathbb{N}}\mathcal{B}_{m},$ with each $%
\mathcal{B}_{m}$ discrete. Then, as $\mathcal{B}$ is a base for $\{Vy:V\in 
\mathcal{V}_{N}\},$ 
\[
X_{N}=\bigcup\nolimits_{m\in \mathbb{N}}\left( \bigcup \{B\in \mathcal{B}%
_{m}:(\exists V\in \mathcal{V}_{N})B\subseteq Vy\}\right) . 
\]%
So for some $m,$ say for $m=M,$ 
\[
\bigcup \{B\in \mathcal{B}_{m}:(\exists V\in \mathcal{V}_{N})B\subseteq Vy\} 
\]%
is non-meagre. But as $\mathcal{B}_{M}$ is discrete, by Banach's Category
Theorem (cf. Prop. B1), there are $\hat{B}\in \mathcal{B}_{M}$ and $\hat{V}%
\in \mathcal{V}_{N}$ with $\hat{B}\subseteq \hat{V}y$ such that $\hat{B}\ $%
is non-meagre. As $\mathcal{V}$ refines $\mathcal{W}$, there is some $\hat{g}%
\in G$ with $\hat{V}\subseteq \hat{g}W,$ so $\hat{B}\subseteq \hat{V}%
y\subseteq \hat{g}Wy,$ and so $\hat{g}Wy$ is non-meagre. As $\hat{g}^{-1}$
is a homeomorphism of $X$, $Wy=Ux$ is also non-meagre in $X$.

If $G$ is analytic, then as $U$ is open, it is also analytic (since open
sets are $\mathcal{F}_{\sigma }$ and Souslin-$\mathcal{F}$ subsets of
analytic sets are analytic, cf. \cite{JayR}), and hence so is $\varphi
_{x}(U).$ Indeed, since $\varphi _{x}$ is continuous and base-$\sigma $%
-discrete, $Ax$ is analytic (Hansell's Theorem, \S 2.1), so Souslin-$%
\mathcal{F}$, and so Baire by Nikodym's Theorem (\S 2.1). $\square $

\bigskip

\textbf{Definition. (}Ancel\textbf{\ }\cite{Anc}\textbf{). }Call the map $%
\varphi _{x}$ \textit{countably-covered}\ if there exist self-homeomorphisms 
$h_{n}^{x}$ of $X$ for $n\in \mathbb{N}$ such that for any open nhd $U$ in $%
G $ the sets $\{h_{n}^{x}(\varphi _{x}(U)):n\in \mathbb{N}\}$ cover $X.$

\bigskip

\textbf{Proposition 1}$^{\prime }$ (cf. Ancell \cite{Anc}) \textit{For the
action }$\varphi :G\times X\rightarrow X$\textit{\ with }$X$\textit{\
non-meagre, if each map }$\varphi _{x}$\textit{\ is countably-covered and
takes open sets to sets with the Baire property, then the action has the
Nikodym property.}

\bigskip

\textbf{Proof.} If $\varphi _{x}$ is countably-covered, then there exist
self-homeomorphisms $h_{n}^{x}$ of $X$ for $n\in \mathbb{N}$ such that for
any open nhd $U$ in $G$ the sets $\{h_{n}^{x}(\varphi _{x}(U)):n\in \mathbb{N%
}\}$ cover $X.$ Then for $X$ non-meagre, there is $n\in \mathbb{N}$ with $%
h_{n}^{x}(\varphi _{x}(U))$ non-meagre, so $Ux=\varphi _{x}(U)$ is itself
non-meagre, being a homeomorphic copy of $h_{n}^{x}(\varphi _{x}(U)).$ As $%
Ux $ is assumed Baire, the action has the Nikodym property. $\square $

\bigskip

For $E$ separable, an immediate consequence of \textit{continuous} maps
taking open sets to analytic sets (which are Baire sets) and of Prop. 1$%
^{\prime }$ is that $\varphi ^{L}$ is a Nikodym action.

For the general context, one needs \textit{demi-open} continuous maps, which
preserve almost completeness (absolute $\mathcal{G}_{\delta }$ sets modulo
meagre sets -- see \cite{Mic2} and its antecedent \cite{Nol}), as it is not
known which linear maps are base-$\sigma $-discrete -- a delicate matter to
determine, since the former include continuous linear surjections (by Lemma
1 below) and preserve almost analyticity as opposed to analyticity.

For present purposes, however, the \textit{monotonicity property} below
suffices. We omit the proof of the following observation (for which see the
opening step in \cite[2.11]{Rud}, or \cite[Ch. 3 \S 12.3]{Con}, or the
Appendix). For the underlying translation-invariant metric of a Fr\'{e}chet
space denote below by $B(a,r)$ the open $r$-ball with centre $a$.

\bigskip

\textbf{Lemma 1.} \textit{For a continuous linear map }$L:X\rightarrow Y$%
\textit{\ from a Fr\'{e}chet space }$X$\textit{\ to a normed space }$Y$%
\textit{, for }$s<t<r$ 
\[
\mathrm{int}(\mathrm{cl}L(B(0,s)))\subseteq L(B(0,t))\subseteq L(B(0,r)). 
\]%
\textit{Hence for }$L(a,r)$\textit{\ convex, either }$L(B(a,r))$\textit{\ is
meagre or differs from }$\mathrm{int}L(B(a,r))$\textit{\ by a meagre set.}

\bigskip

\textbf{Proposition 2.} \textit{For }$L$\textit{\ a continuous linear
surjection\ from a Fr\'{e}chet space }$E$\textit{\ to a non-meagre normed
space }$F,$\textit{\ the action }$\varphi ^{L}$\textit{\ has the Nikodym
property.}

\bigskip

\textbf{Proof. }As in Prop 1$^{\prime }$ for $L:E\rightarrow F$ a continuous
linear surjection, $\{\varphi _{x}^{L}:x\in F\}$ are countably-covered.
Indeed, fixing $x\in F$%
\[
h_{n}^{x}(z):=n(z-x)\text{\qquad (}n\in \mathbb{N}\text{ and }z\in F) 
\]%
is on the one hand a self-homeomorphism satisfying $h_{n}^{x}(\varphi
_{x}(L(V)))=L(nV)$, since $n[(L(v)+x)-x]=nL(v)=L(nv);$ on the other hand the
family%
\[
\{h_{n}^{x}(L(V)+x):n\geq 1\} 
\]%
covers $F,$ as $\{nV:n\in \mathbb{N}\}$ covers $E$ for $V$ any open nhd of
the origin in $E$ (by the `absorbing' property, cf. \cite[4.1.13]{Con}, \cite%
[1.33]{Rud}). In particular, $nL(B(0,1))$ is non-meagre for some $n,$ and so 
$L(B(0,s))$ is non-meagre for any $s$. By Lemma 1, $L(B(0,t))$ for any $t>s$
contains the non-meagre Baire set $\mathrm{cl}L(B(0,s)).$ $\square $

\bigskip

Corollary 1 below is now immediate; it is used in \cite[Th. 5.1]{Ost2} to
prove the `Semi-Completeness Theorem', an Ellis-type theorem \cite[Cor. 2]%
{Ell} (cf. \cite{Ost6}) giving a one-sided continuity condition which
implies that a right-topological group generated by a right-invariant metric
is a topological group.

\bigskip

\textbf{Corollary 1 }(cf. \cite[Th. 5.1]{Ost2}, `Open Homomorphism
Theorem'). \textit{If the continuous surjective homomorphism }$\lambda $%
\textit{\ between normed groups }$G$ \textit{and }$H,$\textit{\ with }$G$%
\textit{\ analytic and }$H$\textit{\ a Baire space, is base-}$\sigma $%
\textit{-discrete, then }$\lambda $\textit{\ is open; in particular, for }$%
\lambda $ \textit{bijective, }$\lambda ^{-1}$\textit{\ is continuous.}

\bigskip

\textbf{Corollary 2}. \textit{For }$L:E\rightarrow F$\textit{\ a continuous
surjective linear map between Fr\'{e}chet spaces, the point evaluations }$%
\varphi _{b}^{L}$\textit{\ for }$b\in F$ \textit{are open, and so }$L$%
\textit{\ is an open mapping.}

\bigskip

\textbf{Proof.} By surjectivity of $L,$ the action is transitive, and by
Prop 2 the action $\varphi ^{L}$ has the Nikodym property. So by Theorem E
above the point-evaluations maps $\varphi _{b}^{L}$ are open. Hence so also
is $L$. $\square $

\section{Proofs}

\subsection{Proof that E $\Longleftrightarrow $ C}

In \cite{BinO1} Th. 3.15 we showed that if the Effros property holds for the
action of a group $G$ on $X$, then $X$ has the crimping property w.r.t. $G$.
We recall the argument, as it is short. Suppose that $x=\lim x_{n}.$ For
each $n,$ take $U=B_{1/n}^{G}(e_{G});$ then $Ux:=\{u(x):u\in U\}$ is an open
nhd of $x,$ and so there exists $h_{n,m}\in U$ with $h_{n,m}(x)=x_{m}$ for
all $m$ large enough, say for all $m>m(n).$ Without loss of generality we
may assume that $m(1)<m(2)<....$ . Put $h_{m}:=e_{G}$ for $m<m(1),$ and for $%
m(k)\leq m<m(k+1)$ take $h_{m}:=h_{k,m}.$ Then $h_{m}\in B_{1/k}^{G}(e_{G}),$
so $h_{m}$ converges to $e_{G}$ and $h_{m}(e_{G})=x_{m}.$

For the converse, suppose that the Effros property fails for $G$ acting on $%
X.$ Then for some open nhd $U$ of $e_{G}$ and some $x\in X,$ $%
Ux:=\{u(x):u\in U\}$ is not an open nhd of $x.$ So for each $n$ there is a
point $x_{n}\in B_{1/n}(x)\backslash Ux.$ As $x_{n}$ converges to $x$ there
are homeomorphisms $h_{n}$ converging to the identity $e_{G}$ with $%
h_{n}(x)=x_{n}$. As $U$ is an open nhd of $e_{G}$ and since $h_{n}$
converges to $e_{G}$, there is $N$ such that $h_{n}\in U$ for $n>N.$ In
particular, for any $n>N,$ $h_{n}(x)=x_{n}\in Ux,$ a contradiction.

\subsection{Weak S}

We view Th. S as having `two tasks': to find a `translator of the sequence' $%
\tau ,$ and to locate it in a given Baire non-meagre subset of the group --
provided that subset satisfies a consistency condition (a necessary
condition).

For clarity we break the tasks the into two steps -- the first delivering a
weaker version of S in Proposition 3 below. The arguments are based on the
following lemma. We note a corollary, observed earlier by van Mill in the
case of metric topological groups (\cite[Prop. 3.4]{vMil2}), which concerns
a co-meagre set, but we need its refinement to a localized version for a
non-meagre set.

\bigskip

\textbf{Separation Lemma.} \textit{Let }$G$\textit{\ be a normed group, with
separately continuous and transitive Nikodym action on a non-meagre space }$%
X $\textit{. Then for any point }$x$\textit{\ and any }$F$\textit{\ closed
nowhere dense, }$W_{x,F}:=\{\alpha \in G:\alpha (x)\notin F\}$\textit{\ is
dense open in }$G$\textit{. In particular, }$G$ \textit{separates points
from nowhere dense closed sets.}

\bigskip

\textbf{Proof.} The set $W_{x,F}$ is open, being of the form $\varphi
_{x}^{-1}(X\backslash F)$ with $\varphi _{x}$ continuous (by assumption). By
the Nikodym property, for $U$ any non-empty open set in $G,$ the set $Ux$ is
non-meagre, and so $Ux\backslash F$ is non-empty, as $F$ is meagre. But then
for some $u\in U$ we have $u(x)\notin F$. $\square $

\bigskip

\textbf{Corollary 2.} \textit{If }$G$\textit{\ is a normed group, Baire in
the norm topology with transitive and separately continuously Nikodym action
on a non-meagre space }$X$ \textit{space, and }$T$ \textit{is co-meagre in }$%
X$\textit{-- then for countable }$D\subseteq X,$\textit{\ the set }$%
\{g:g(D)\subseteq T\}$\textit{\ is a dense }$\mathcal{G}_{\delta }$\textit{.}

\textit{In particular, this holds if }$G$\textit{\ is analytic and each
point-evaluation map }$\varphi _{x}:g\rightarrow g(x)$\textit{\ is base-}$%
\sigma $\textit{-discrete.}

\bigskip

\textbf{Proof. }Without loss of generality, the co-meagre set is of the form 
$T=U\backslash \bigcup_{n\in \omega }F_{n}$ with each $F_{n}$ closed and
nowhere dense, and $U$ open. Then, by the Separation Lemma and as $G$ is
Baire, 
\[
\{g\in G:g(D)\subseteq T\}=\bigcap\nolimits_{n\in \omega }\{g:g(D)\cap
F_{n}=\emptyset \}=\bigcap\nolimits_{d\in D,n\in \omega }\{g:g(d)\notin
F_{n}\} 
\]%
is a dense $\mathcal{G}_{\delta }$. $\square $

\bigskip

\textbf{Proposition 3.} \textit{If }$T$\textit{\ is a Baire non-meagre
subset of a metric space }$X$\textit{\ and }$G$\textit{\ a normed group,
Baire in its norm topology, acting separately continuously and transitively
on }$X,$ \textit{with the Nikodym property -- then, for every convergent
sequence }$x_{n}$\textit{\ with limit }$x_{0}$\textit{\ there is }$\tau \in
G $\textit{\ and an integer }$N$\textit{\ with }$\tau x_{0}\in T$\textit{\
and}%
\[
\{\tau (x_{n}):n>N\}\subseteq T. 
\]

\textbf{Proof.} Write $T:=M\cup (U\backslash \bigcup\nolimits_{n\in \omega
}F_{n})$ with $U$ open, $M$ meagre and each $F_{n}$ closed and nowhere dense
in $X$. Let $u_{0}\in T\cap U.$ By transitivity there is $\sigma \in G$ with 
$\sigma x_{0}=u_{0}.$ Put $u_{n}:=\sigma x_{n}.$ Then $u_{n}\rightarrow
u_{0}.$ Put%
\[
C:=\bigcap\nolimits_{m,n\in \omega }\{\alpha \in G:\alpha (u_{m})\notin
F_{n}\}, 
\]%
a dense $\mathcal{G}_{\delta }$ in $G;$ then, by the Separation Lemma above,
as $G$ is Baire,%
\[
\{\alpha \in G:\alpha (u_{0})\in U\}\cap C 
\]%
is non-empty. For $\alpha $ in this set we have $\alpha (u_{0})\in
U\backslash \bigcup\nolimits_{n\in \omega }F_{n}.$ Now $\alpha
(u_{n})\rightarrow \alpha (u_{0}),$ by continuity of $\alpha $, and $U$ is
open. So for some $N$ we have for $n>N$ that $\alpha (u_{n})\in U.$ Since $%
\{\alpha (u_{m}):m=1,2,..\}\in X\backslash \bigcup\nolimits_{n\in \omega
}F_{n},$ we have for $n>N$ that $\alpha (u_{n})\in U\backslash
\bigcup\nolimits_{n\in \omega }F_{n}\subseteq T.$

Finally put $\tau :=\alpha \sigma ;$ then $\tau (x_{0})=\alpha \sigma
(x_{0})\in T$ and $\{\tau (x_{n}):n>N\}\subseteq T.$ $\square $

\subsection{Proof of S}

We work in the right norm topology and use the notation of the preceding
proof (of Proposition 3), so that $U$ here is the quasi-interior of $T$ and $%
\sigma x_{0}=u_{0}.$ As $e_{G}\in A^{q}$ and $A$ is a non-meagre Baire set,
we may without loss of generality write $A=B_{\varepsilon }(e_{G})\backslash
\bigcup\nolimits_{n}G_{n},$ where each $G_{n}$ is closed nowhere dense with $%
e_{G}\notin G_{n}$ and $B_{\varepsilon }(e_{G})$ is the quasi-interior of $%
A. $

As $A^{q}x_{0}\cap T^{q}$ is non-empty, there is $\alpha _{0}\in
B_{\varepsilon }(e_{G})$ with $\alpha _{0}x_{0}\in U$ (but, we want a better 
$\alpha $ so that $\alpha x_{0}\in T$ and $\alpha \in A).$ Put $\beta
_{0}=\alpha _{0}\sigma ^{-1};$ then 
\begin{eqnarray*}
\beta _{0} &=&\alpha _{0}\sigma ^{-1}\in B_{\varepsilon }(e_{G})\sigma
^{-1}\cap \{\alpha :\alpha (x_{0})\in U\}\sigma ^{-1} \\
&=&B_{\varepsilon }(e_{G})\sigma ^{-1}\cap \{\beta :\beta (\sigma x_{0})\in
U\}=B_{\varepsilon }(e_{G})\sigma ^{-1}\cap \{\beta :\beta (u_{0})\in U\},
\end{eqnarray*}%
i.e. the open set $\{\beta :\beta (u_{0})\in U\}\cap B_{\varepsilon
}(e_{G})\sigma ^{-1}$ is non-empty. So%
\[
(C\backslash \bigcup\nolimits_{n}G_{n}\sigma ^{-1})\cap \{\beta :\beta
(u_{0})\in U\}\cap B_{\varepsilon }(e_{G})\sigma ^{-1}\neq \emptyset , 
\]%
since $G$ is a Baire space and each $G_{n}\sigma ^{-1}$ is closed and
nowhere dense in $G$ (as the right shift $g\rightarrow g\sigma ^{-1}$ is a
homeomorphism).

So there is $\beta $ with $\beta (u_{0})\in U$ such that $\alpha :=\beta
\sigma \in B_{\varepsilon }(e_{G})\backslash \bigcup\nolimits_{n}G_{n}=A.$
That is, $\alpha x_{0}=\beta u_{0}\in U;$ so $\beta (u_{n})\in U$ for large $%
n$, for $n>N\ $say, as $\alpha x_{0}=\lim \alpha x_{n}=\lim \beta \sigma
x_{n}=\lim \beta u_{n}$. But $\{\beta (u_{m}):m=1,2,..\}\in X\backslash
\bigcup\nolimits_{n}F_{n},$ as $\beta \in C;$ so $\beta (u_{n})\in
U\backslash \bigcup\nolimits_{n}F_{n}\subseteq T$ for $n>N$.

Finally, $\alpha (x_{0})=\beta \sigma (x_{0})\in T$ and $\{\alpha
(x_{n}):n>N\}\subseteq T.$ $\square $

\subsection{\textbf{P}roof that S $\Longrightarrow $ E}

Assume $G$ acts transitively on $X$ and that $X$ is non-meagre. Let $%
B:=B_{\varepsilon }(e_{G})$ and suppose that for some $x$ the set $Bx$ is
not a nhd of $x.$ Then there is $x_{n}\rightarrow x$ with $x_{n}\notin Bx$
for each $n.$ Take $A:=B_{\varepsilon /2}(e_{G})$ and note first that $A$ is
a symmetric open set ($A^{-1}=A,$ since $||g||=||g^{-1}||$), and secondly
that by the Nikodym property $Ax$ contains a non-meagre, Baire subset $T$.
So by Theorem S, as $Ax$ meets $T^{q},$ there are $a\in A$ (which being open
has the Baire property) and a co-finite $\mathbb{M}_{a}$ such that $%
ax_{m}\in Ax$ for $m\in \mathbb{M}_{a}$. For any such $m,$ choose $b_{m}\in
A $ with $ax_{m}=b_{m}x.$ Then $x_{m}=a^{-1}b_{m}x\in A^{2}x\subseteq Bx,$ a
contradiction (note that $a^{-1}\in A,$ by symmetry).

As earlier, in the special case that $G$ is (metrizable and) analytic, $A$
is analytic, since open sets are $\mathcal{F}_{\sigma }$ and Souslin-$%
\mathcal{F}$ subsets of analytic sets are analytic, cf. \cite[Th. 2.5.3]%
{JayR}, by Prop. 3 $Ax$ is Baire non-meagre, as $\varphi _{x}$ is base-$%
\sigma $-discrete.

\section{Complements}

\noindent \textit{1. Normed groups.} For $T$ an \textit{algebraic} group
with neutral element $e=e_{G}$, say that $||\cdot ||:T\rightarrow \mathbb{R}%
_{+}$ is a \textit{group}-\textit{norm }(\cite{ArhT}, \cite{BinO1}) if the
following properties hold:

(i) \textit{Subadditivity} (Triangle inequality): $||st||\leq ||s||+||t||;$

(ii) \textit{Positivity}: $||t||>0$ for $t\neq e$ and $||e||=0;$

(iii) \textit{Inversion} (Symmetry): $||t^{-1}||=||t||.$

\noindent Then $(T,||.||)$ is called a \textit{normed-group.}

The group-norm generates a right and a left \textit{norm topology} via the
right-invariant and left-invariant metrics $d_{R}^{T}(s,t):=||st^{-1}||$ and 
$d_{L}^{T}(s,t):=||s^{-1}t||=d_{R}^{T}(s^{-1},t^{-1}).$ In the right
norm-topology the right shift $\rho _{t}(s):=st$ is a uniformly continuous
homeomorphism, since $d_{R}(sy,ty)=d_{R}(s,t)$, so in particular the group
is a right topological group; likewise in the left norm-topology the left
shift $\lambda _{s}(t)=st$ is a uniformly continuous homeomorphism. Since $%
d_{L}^{T}(t,e)=d_{L}^{T}(e,t^{-1})=d_{R}^{T}(e,t),$ convergence at $e$ is
identical under either topology. One may refer to whichever topology is
appropriate, since despite their differences, they are homeomorphic. Thus
one may say, for instance, that separability is a \textit{norm property}, as
separability of either norm topology implies that of the other.

See \cite{Ost5} for a characterization of almost-complete normed groups as
non-meagre almost-analytic (i.e. analytic modulo a meagre set).

\bigskip 

\noindent \textit{2. Genericity }Our focus on various shift-compactness
theorems derives from their close affinity with the literature of `generic'
automorphisms (for which see \cite{Mac}, also described in the introduction
to  \cite{Ost5}) and from a multitude of its applications, for which see
e.g. \cite{BinO1}, and these offering a unifying sequential
compactness-like, or combinatorial, perspective on category-measure duality
and on other, apparently unrelated, problems.

\bigskip 

\noindent \textit{3. Assumptions in Theorem S. }With regard to the
assumption of separate continuity, note that a theorem of Bouziad (\cite[Th.
3]{Bou}) implies that a separately continuous action by a metrizable
left-topological Baire group acting on a metric space is in fact jointly
continuous; we retain the only apparently weaker hypothesis of separate
continuity, because there are variants of Theorem S, where joint continuity
is absent (and the group is not Baire): see \cite{MilO}. The result is
connected with van Mill's \textit{Separation Property}: say that $X$ has SP (%
\cite{vMil2}) with respect to a group of homeomorphisms $\mathcal{G}$ if for
any countable set $D$ in $X$ and any meagre set $M$ in $X,$ there is a
homeomorphism $g\in \mathcal{G}$ such that $g(D)\cap M=\emptyset ;$ it is
interesting to recall here the Galvin-Mycielski-Solovay characterization of
subsets of $\mathbb{R}$ of strong measure zero as sets $A$ such that for any
meagre $M$ there is $g$ with $(g+A)\cap M=\emptyset $ (see \cite[\S 3]{Mil}%
). Stated equivalently, SP asserts that for $D$ countable and $T$ co-meagre,
there is $g\in \mathcal{G}$ with $g(D)\subseteq T$ (see also the next
remark). Compared to this restatement, Theorem S refers on the one hand to a
smaller class of countable sets (convergent sequences, or their co-finite
parts), but on the other hand asserts embeddability into a larger class of
sets -- sets $T$ that are `locally' rather than globally co-meagre; for
further information see also \cite{MilO}.

\bigskip

\noindent \textit{4. Shift-compactness and the SP property. }In view of the
strength of Th. S and to place our results in context, we briefly summarize
some of the relevant results of \cite[Th. 1.1]{vMil2}. A (separable metric)
space with SP is Baire, by the proof which we imported for Propositions B1
and B2 in \S 2; likewise, an almost complete non-meagre space with SP is
completely Baire. Since an absolutely co-analytic space is Polish iff it is
`completely Baire', i.e. closed-hereditarily Baire (closed subspaces are
Baire), for which see Kechris \cite[Cor. 21.21]{Kec2}, it follows that an
absolutely Borel space with the SP is Polish (\cite[Th. 1.1]{vMil2}). More
generally, if an analytic group $G$ acts on a space $X$ and SP holds w.r.t. $%
G$, then $X$ is Polish. Van Mill also shows from his Prop 3.4 (cf. our Prop.
2) that a locally compact homogeneous space has the SP. It seems likely
that, just as with Proposition B, more of these arguments can be copied
across in the language of Th. S.

It is noted in \cite{BinO2} that a Polish space which is strongly locally
homogeneous has SP.

\bigskip

\noindent 5. \textit{Index-}$\sigma $\textit{-discrete maps. }In many
circumstances it is easier to work with index-$\sigma $-discrete maps: see 
\cite{Ost1} for a brief discussion of this point and for relations with the
automatic continuity of homomorphisms (noted earlier in \cite{Nol}).
Furthermore, if $f:X\rightarrow Y$ is injective and \textit{closed-analytic,}
i.e. carries closed sets to analytic sets (or, alternatively \textit{%
open-analytic}, mutatis mutandis), then $f$ is base-$\sigma $-discrete -- in
fact index-$\sigma $-discrete (\cite{Han3} Prop. 3.14). On the other hand,
if $f$ is surjective and closed with $Y$ metrizable, then $f$ is base-$%
\sigma $-discrete (\cite{Han3} Prop. 3.10). So if the group action is such
that each $\varphi _{x}$ is open-analytic and \textit{injective} (as when a
vector space acts on itself), then each $\varphi _{x}$ is index- and so base-%
$\sigma $-discrete.

\bigskip

\noindent 6. \textit{Generalizations of analyticity. }A Hausdorff space is 
\textit{almost analytic} if it is analytic modulo a meagre set. Similarly, a
space $X^{\prime }$ is absolutely $\mathcal{G}_{\delta }$, or an \textit{%
absolute-}$\mathcal{G}_{\delta }$, if $X^{\prime }$ is a $\mathcal{G}%
_{\delta }$ in all spaces $X$ containing $X^{\prime }$ as a subspace. This
latter property is equivalent to complete metrizability in the realm of
metrizable spaces \cite[Th. 4.3.24]{Eng} (and to topological/\v{C}ech
completeness in the realm of completely regular spaces -- \cite[\S 3.9]{Eng}%
). So a metrizable absolute-$\mathcal{G}_{\delta }$ is analytic. A metric
space is \textit{almost complete} if it contains a dense absolute $\mathcal{G%
}_{\delta }$. The notion of `almost completeness' is due to Frol\'{\i}k in 
\cite{Frol} (but its name to Michael \cite{Mic2} -- see also \cite{AaL} and 
\cite{BinO1}).

Consequently, it is natural to assume that when a group $G$ acts on a space $%
X$ each point evaluation map $\varphi _{x}:g\rightarrow g(x)$ is not only
continuous but also base-$\sigma $-discrete; the \textit{Nikodym property}
above is implied by the former property when $G\ $is analytic. Note that 
\cite{Han3} Ex. 3.12 shows that, for $D=(\mathbb{R},d_{\text{discrete}})$
the discrete space of cardinality the continuum, the projection from $%
D\times \lbrack 0,1]\rightarrow \lbrack 0,1]$ is an open mapping that is not
base-$\sigma $-discrete (by reference to the closed discrete graph of a
bijection between $D$ and $[0,1])$.

Index-$\sigma $-discrete maps (easier to work with) are base-$\sigma $%
-discrete (\cite{Han3} Prop. 3.7(i)); the latter combined with continuity
preserves analyticity, as mentioned earlier. Evidently, if all evaluation
maps $\varphi _{x}$ are continuous and base-$\sigma $-discrete and $A$ is an
open subset of an analytic group $G,$ acting on $X,$ then $A$ is analytic,
and so $\varphi _{x}(A)=Ax$ is analytic; thus the point evaluations are then
open-analytic. So our Nikodym hypothesis on the action is not far from
demanding that point evaluations be open-analytic.

\bigskip 

\noindent 7. \textit{Convexity and Baire Property. }We have seen above that
analyticity confers the Baire property. Convexity may also confer it: for $A$
convex, if \textrm{int}$A,$ its interior, is non-empty, then \textrm{int}$A$
is dense in $A,$ and so $A$ has the Baire property: $A\backslash \mathrm{int}%
(A)\subseteq \mathrm{cl}(A)\backslash \mathrm{int}(\mathrm{cl}(A)),$ a
nowhere dense set; so $A$ differs from its interior by a nowhere dense
subset of $\mathrm{cl}(A)\backslash \mathrm{int}(\mathrm{cl}(A)).$ If $A$ is
non-meagre its closure has non-empty interior, raising the question of
whether $A$ itself has non-empty interior.

The Banach context is particularly transparent; Banach's lemma (see
Appendix) yields: \textit{for a continuous linear map }$T:X\rightarrow Y$%
\textit{\ from a Banach space }$X$\textit{\ to a normed space }$Y$\textit{\
and }$x\in X,$\textit{\ if }$B(0,\rho )\subseteq $\textrm{cl}$(T(\bar{B}%
(x,r)))$\textit{\ with }$r,\rho >0$\textit{, then} $B(0,\rho )\subseteq T(%
\bar{B}(x,r)).$

So, as above, by translation, interpret this lemma as asserting that for $r>0
$ and $A:=L(\bar{B}(0,r))$ one has $\mathrm{int}(\mathrm{cl}(A))\subseteq 
\mathrm{int}(A).$ So $A\backslash \mathrm{int}(A)\subseteq \mathrm{cl}%
(A)\backslash \mathrm{int}(\mathrm{cl}(A)),$ a nowhere dense set, so $A$
differs from its interior by a nowhere dense subset of $\mathrm{cl}%
(A)\backslash \mathrm{int}(\mathrm{cl}(A)),$ and so has the Baire property.

In the Fr\'{e}chet case, with the balls referring to the underlying
translation-invariant metric, one may apply Lemma 1 of \S 2.3. Taking the
convex set $B(0,r)$ with $r>0$, suppose $T(B(0,r))$ is non-meagre; then $%
T(B(0,s))$ is non-meagre for some $0<s<r.$ By the claim, for $s<t<r,$%
\[
\emptyset \neq \mathrm{int}(\mathrm{cl}T(B(0,s)))\subseteq T(B(0,t)).
\]%
As $T(B(0,t))$ is convex, $\mathrm{int}T(B(0,t))$ is convex, and being
non-empty is dense in $\mathrm{cl}T(B(0,t)),$ i.e. $\mathrm{cl}T(B(0,t))=%
\mathrm{cl}(\mathrm{int}T(B(0,t))).$ So $T(B(0,t))$ differs from $\mathrm{int%
}T(B(0,t))$ by a meagre set, as%
\begin{eqnarray*}
T(B(0,t))\backslash \mathrm{int}T(B(0,t)) &\subseteq &\mathrm{cl}%
T(B(0,t))\backslash \mathrm{int}T(B(0,t)) \\
&=&\mathrm{cl}(\mathrm{int}T(B(0,t)))\backslash \mathrm{int}T(B(0,t)).
\end{eqnarray*}%
Finally, taking $s<t_{n}\nearrow r,$ one has 
\[
T(B(0,r))\backslash \mathrm{int}T(B(0,r))\subseteq \bigcup\nolimits_{n\in 
\mathbb{N}}\mathrm{cl}(\mathrm{int}T(B(0,t_{n})))\backslash \mathrm{int}%
T(B(0,t_{n})).
\]%
So for $B(0,r)$ convex,\textit{\ either }$T(B(0,r))$\textit{\ is meagre or
differs from }$\mathrm{int}T(B(0,r))$\textit{\ by a meagre set.}

\bigskip 

\noindent 8. \textit{Topological and semitopological groups. }When a group $G
$ equipped with a topology such that the (action) map $(g_{1},g_{2})%
\rightarrow g_{1}g_{2}$ from $G^{2}$ to $G$ is separately continuous (i.e.
left and right shifts are continous), then it is said to be a \textit{%
semitopological group}. Thus a topological group is semitopological. If also
the topology is metrizable and $G$ is Baire, Theorem S applies and asserts
here that if $T$ is non-meagre in $G$ and $g_{n}\rightarrow g,$ then for
some $t\in T$ the point $tg$ is in $T$ and almost all of the sequence $tg_{n}
$ is in $T.$ As in Theorem S, under these additional assumptions, group
multiplication in $G$ is jointly continuous; indeed, Bouziad \cite{Bou}
proves that a semitopological Baire $p$-space (and metric spaces are $p$%
-spaces) has jointly continuous multiplication (is `paratopological').

\bigskip

\noindent \textit{9. Specializing the proof of Theorem S to semitopological
groups.} The argument proving Th. S is particularly transparent in the case
of a semitopological group and follows a standard real-line argument, as
follows. Suppose $x_{n}\rightarrow x_{0}$ and consider $T:=U\backslash
\bigcup F_{m}$ with $F_{m}$ closed and nowhere dense, and $U$ open. For any $%
u_{0}\in T\cap U,$ take $\sigma (x)=xx_{0}^{-1}u_{0},$ which is a continuous
right-shift, so that $u_{n}=\sigma (x_{n})=x_{n}(x_{0}^{-1}u_{0})\rightarrow
u_{0}=\sigma (x_{0}).$ Now%
\[
u_{0}^{-1}U\cap \bigcap\nolimits_{n,m}u_{n}^{-1}(X\backslash F_{m})
\]%
is a dense $\mathcal{G}_{\delta }$, since left-shifts are homeomorphisms. As 
$G$ is Baire, there is $g$ with 
\[
g\in u_{0}^{-1}U\cap \bigcap\nolimits_{n,m}u_{n}^{-1}(X\backslash F_{m}).
\]%
Then $u_{0}g\in U$ and $u_{0}g\notin F_{m}$ for all $m,$ so $u_{0}g\in T.$
By continuity of right-shifts, $u_{n}g\rightarrow u_{0}g$, so for large $n,$
say for $n>N,$ we have%
\[
u_{n}g\in U.
\]%
For each such $n$ we have $u_{n}g\notin F_{m}$ for all $m$ and so $u_{n}g\in
T.$ Thus, as $u_{n}g=x_{n}(x_{0}^{-1}u_{0})g,$ we conclude the existence for
some $g\in G$ and $u_{0}\in T$ with $t:=u_{0}g\in T$ that, as in Theorem S%
\[
x_{n}x_{0}^{-1}t\in T\text{ for }n>N.
\]

\noindent \textit{10. From E\ back to S.} In Proposition 3 density of the
set $W_{x,F}$ was deduced from a Lemma which, as noted, asserts a very weak
microtransitivity. Unsurprisingly, Property E also implies density of that
set (see below). However, Property E does not imply that $G$ is Baire, since
(see Introduction) there do exist meagre analytic groups acting on a
non-meagre space which by van Mill's result have Property $E$
notwithstanding. So we can go no further with the density argument of
Proposition 3. Of course, as in Proposition 3, if we know that $G$ is Baire,
then E implies S.

As for the density claim, consider $\beta \in G.$ Suppose first that $\beta
(x)\in F.$ By the Effros property, the set $B_{\varepsilon }(\beta
)x=B_{\varepsilon }(e_{G})\beta x$ is an open nhd of $\beta x.$ As $F$ is
nowhere dense and closed, there is $y\in $ \textrm{int}$(B_{\varepsilon
}(e_{G})\beta x)\backslash F$. So there is $\gamma \in B_{\varepsilon
}(e_{G})$ such that $y=\alpha x:=\gamma \beta x\notin F.$ So $a=\gamma \beta
\in B_{\varepsilon }(\beta )\cap \{\alpha :\alpha (x)\notin F\}.$

If, on the other hand, $\beta (x)\notin F,$ then $\beta \in B_{\varepsilon
}(\beta )\cap \{\alpha :\alpha (x)\notin F\}.$ This proves the claim.

\bigskip 

\textbf{A personal note. }Whilst the present author's entry into mathematics
owes hugely both to Karol Borsuk and Ambrose Rogers (thesis advisor), being
confirmed as a topologist is down to a first meeting with Mary Ellen Rudin
in 1972 (at the Keszthely conference) and subsequent frequent stays at UW
Madison, visiting her and the wonderous set-theoretic community there. It is
thus a pleasure to dedicate this paper especially to her memory and to
express once more the great debt to UW friends and colleagues, among whom
also was Anatole Beck, recently passed away.

\bigskip 

\textbf{Acknowledgement. }I thank the Referee for some thought-provoking
comments that influenced the final presentation, and Henryk Toru\'{n}czyk
for drawing my attention to Ancel's work and related literature.

\bigskip 

\noindent Mathematics Department, London School of Economics, Houghton
Street, London WC2A 2AE \newline
a.j.ostaszewski@lse.ac.uk\newpage 

\section*{Appendix}

\noindent For convenience we include here the

\bigskip 

\noindent \textbf{Lemma.} (see e.g. \cite[Lemma 5.4]{TayL}; \cite[Ch. III;
Proof of Th.3]{Ban}) \textit{For a continuous linear map }$T:X\rightarrow Y$%
\textit{\ from a Banach space }$X$\textit{\ to a normed space }$Y$\textit{,
if }$B(0,\rho )\subseteq $\textrm{cl}$(T(\bar{B}(a,r))),$\textit{\ then} $%
B(0,\rho )\subseteq T(\bar{B}(a,r)).$

\bigskip 

\noindent \textbf{Proof.} Let $\varepsilon $ be arbitrary with $%
0<\varepsilon <1$ and suppose that $B(0,\rho )\subseteq $\textrm{cl}$(T(\bar{%
B}(a,r)))$. If $y\in B(0,\rho ),$ then $||y-y_{1}||<\varepsilon \rho $ for
some $y_{1}=Tx_{1}$ with $x_{1}\in B(a,r)$. So $x_{1}=a+z_{1}$ with $%
||z_{1}||<r.$ But $||y-Tx_{1}||<\varepsilon \rho $ and, by homogeneity of
the norm in $Y$,%
\[
B(0,\varepsilon \rho )=\varepsilon B(0,\rho )\subseteq \mathrm{cl}T(\bar{B}%
(\varepsilon a,\varepsilon r)),
\]%
so 
\[
||y-Tx_{1}-Tx_{2}||<\varepsilon ^{2}\rho .
\]%
for some $x_{2}\in \bar{B}(\varepsilon a,\varepsilon r),$ i.e. $%
x_{2}=\varepsilon a+z_{2}$ for some $z_{2}$ with $||z_{2}||\leq \varepsilon
r.$ Continue by induction. Now%
\[
x_{1}+...+x_{n}=a(1+\varepsilon +\varepsilon ^{2}+...+\varepsilon
^{n-1})+z_{1}+..+z_{n},
\]%
and%
\[
||z_{1}+...+z_{n}||\leq ||z_{1}||+...+||z_{n}||\leq r(1+...+\varepsilon
^{n-1}).
\]%
By completeness, $x_{1}+...+x_{n}\rightarrow x_{\ast }\in \bar{B}%
(a/(1-\varepsilon ),r/(1-\varepsilon )$ and $Tx_{\ast }=y.$ So 
\[
B(0,\rho )\subseteq T(\bar{B}(a/(1-\varepsilon ),r/(1-\varepsilon ))
\]%
and so%
\[
B(0,(1-\varepsilon )\rho )=(1-\varepsilon )B(0,\rho )\subseteq T(\bar{B}%
(a,r)).
\]%
As $\varepsilon $ was arbitrary we conclude that%
\[
B(0,\rho )=\bigcup\nolimits_{0<\varepsilon <1}B(0,(1-\varepsilon )\rho
)\subseteq T(\bar{B}(a,r)).
\]%
In particular $0\in \mathrm{int}(T(\bar{B}(a,r))).$ $\square $ \newpage 

\end{document}